\documentclass[12pt,a4paper]{article}

\usepackage[latin2]{inputenc}
\usepackage{amssymb}
\usepackage{paralist}
\usepackage{amsmath, calc}
\usepackage{amsfonts}
\usepackage{amsthm}
\usepackage{fullpage}
\usepackage{enumerate}
\usepackage{paralist}

\newtheorem{thm}{Theorem}
\newtheorem{lem}{Lemma}

\newtheorem{conj}{Conjecture}

\newtheorem{rem}{Remark}

\begin{document}

\title{On generalized Stanley sequences}
\author{S\'andor Z. Kiss \thanks{Institute of Mathematics, Budapest
University of Technology and Economics, H-1529 B.O. Box, Hungary;
kisspest@cs.elte.hu;
This author was supported by the National Research, Development and Innovation Office NKFIH Grant No. K115288.}, Csaba
S\'andor \thanks{Institute of Mathematics, Budapest University of
Technology and Economics, H-1529 B.O. Box, Hungary, csandor@math.bme.hu.
This author was supported by the OTKA Grant No. K109789. This paper was supported
by the J\'anos Bolyai Research Scholarship of the Hungarian Academy of Sciences.}, Quan-Hui Yang \thanks{School of Mathematics and Statistics, Nanjing University of Information Science and Technology, Nanjing 210044, China; yangquanhui01@163.com; This author was supported by the National Natural Science
Foundation for Youth of China, Grant No. 11501299, the Natural
Science Foundation of Jiangsu Province, Grant Nos.
BK20150889,~15KJB110014 and the Startup Foundation for Introducing
Talent of NUIST, Grant No. 2014r029.}
}
\date{}
\maketitle

\begin{abstract}
\noindent
Let $\mathbb{N}$ denote the set of all nonnegative integers.
Let $k\ge 3$ be an integer and $A_{0} = \{a_{1}, \dots{}, a_{t}\}$ $(a_{1} < \ldots< a_{t})$ be a nonnegative set
which does not contain an arithmetic progression of length $k$. We denote $A = \{a_{1}, a_{2}, \dots{}\}$
defined by the following
greedy algorithm: if $l \ge t$ and $a_{1}, \dots{}, a_{l}$ have already been defined, then $a_{l+1}$ is the smallest integer $a > a_{l}$ such that
$\{a_{1}, \dots{}, a_{l}\} \cup \{a\}$ also does not contain a $k$-term arithmetic progression.
This sequence $A$ is called the Stanley sequence of order $k$ generated by
$A_{0}$. In this paper, we prove some results about various generalizations of the Stanley sequence.

 {\it
2010 Mathematics Subject Classification:} Primary 11B75.

{\it Keywords and phrases:}  Stanley sequences, arithmetic progression, probabilistic method
\end{abstract}


\section{Introduction}
Let $\mathbb{N}$ denote the set of all nonnegative integers. For a finite set
$A_{0} \subset \mathbb{N}$, $A_{0} = \{a_{1}, \dots{}, a_{t}\}$
$(a_{1} < \dots{} < a_{t})$ which does not contain an arithmetic progression of length $k$, we denote $A = \{a_{1}, a_{2}, \dots{}\}$ the sequence
defined by the following greedy
algorithm: if $l \ge t$ and $a_{1}, \dots{}, a_{l}$ have already been defined, then $a_{l+1}$ is the smallest integer $a > a_{l}$ such that
$\{a_{1}, \dots{}, a_{l}\} \cup \{a\}$ does not contain an arithmetic progression of length $k$. This sequence is called Stanley sequence of
order $k$ generated by $A_{0}$.

\begin{rem}
If $k = 3$, then $A$ is a Stanley sequence of order $3$ if and only if $n \in A~\Leftrightarrow ~n \ne 2b-a$, where $a$, $b < n$ and $a, b \in A$.
\end{rem}

To investigate the density of sets without arithmetic progressions is
one of the most popular topic in additive combinatorics.  In 1953,
Roth \cite{rot} proved that every subset of the set of integers with positive
upper density contains an arithmetic progression of length three. On
the other hand, Behrend \cite{ber} constructed a
dense set without any arithmetic progression of length three. The name
Stanley sequences established by Erd\H{o}s et al. \cite{erd} and the definition
originates with Odlyzko and Stanley from 1978. In their joint paper,
they \cite{ost} constructed sets without arithmetic progression of length
three by using the greedy algorithm. In this paper we generalize the concept of Stanley
sequences in two directions. First, we will define the $AP_{k}$ - covering sequences.
In the first three theorems, we study the density of these sequences. In the other direction,
we extend the definition of Stanley sequence according to Remark 1. In the last theorem we
give a fully description of the structure of such sets when $A_{0} = \{a_{0}\}$.
Now we give the notations and definitions we are working with.

Let $A(n)$ be the number of elements of $A$ up to $n$ i.e.,
\[
A(n) = \sum_{\overset{a \in A}{a \le n}}1.
\]
We denote $f = O(g)$ by $f \ll g$.
Gerver and Ramsey \cite{gra} proved that
if $A$ is a Stanley sequence of order $3$, then
\[
\liminf_{n \rightarrow \infty}\frac{A(n)}{\sqrt{n}} \ge \sqrt{2}.
\]
A few years later, Moy \cite{moy} rediscovered this inequality. Recently Chen and Dai \cite{chd} proved that if $A$ is a Stanley sequence of order $3$,
then
\[
\limsup_{n \rightarrow \infty}\frac{A(n)}{\sqrt{n}} \ge 1.77.
\]

We say a sequence $A \subseteq \mathbb{N}$ is an {\em $AP_{k}$ - covering sequence} if
there exists an integer $n_{0}$ such that if $n > n_{0}$, then there exist $a_{1} \in A,\dots{}, a_{k-1} \in A$, $a_{1} < \cdots < a_{k-1} < n$ such that $a_{1}, \dots{}, a_{k-1}, n$ form a $k$-term arithmetic progression. Clearly, if $A$ is a Stanley sequence
 of order $k$, then $A$ is also an $AP_k$ - covering sequence.

 Using Gerver and Ramsey' idea, we can give a lower bound for $A(n)$ if $A$ is an
$AP_{k}$-covering sequence. Obviously
\begin{eqnarray*}
&&n-n_0\le |\{(a_m,b_m):n_0<m\le n,a_m,b_m<m,a_m,b_m\in A,\\
&&a_m,b_m,m\mbox{ form an arithmetic progression of length three}\}|
\le \binom{A(n)}{2}.
\end{eqnarray*}
Hence we have $A(n)\ge \sqrt{2n-2n_0+0.25}+0.5$, which implies
\[
\liminf_{n \rightarrow \infty}\frac{A(n)}{\sqrt{n}} \ge \sqrt{2}.
\]
Similarly, using Chen and Dai's proof, we may verify that if $A$ is an $AP_{3}$ - covering sequence, then
\[
\limsup_{n \rightarrow \infty}\frac{A(n)}{\sqrt{n}} \ge 1.77.
\]
We omit the details.

In this paper we prove the following theorems.

\begin{thm}
There exists an $AP_{3}$ - covering sequence $A$ such that
\[
\liminf_{n \rightarrow \infty}\frac{A(n)}{\sqrt{n}} \le 2.
\]
\end{thm}

\begin{thm}
There exists an $AP_{3}$ - covering sequence $A$ such that
\[
\limsup_{n \rightarrow \infty}\frac{A(n)}{\sqrt{n}} \le 34.
\]
\end{thm}

\begin{thm}
There exists an $AP_{k}$ - covering sequence $A$ such that
\[
A(n) \ll_{k} (\log n)^{1/(k-1)}n^{\frac{k-2}{k-1}}.
\]
\end{thm}
We pose the following conjecture.

\begin{conj}
\begin{itemize}
\item[(i)] For any integer $k \ge 3$, there exists an $AP_{k}$ - covering sequence $A$ such that
\[
\limsup_{n \rightarrow \infty}\frac{A(n)}{n^{\frac{k-2}{k-1}}} < \infty.
\]
\item[(ii)] For any $AP_{k}$ - covering sequence $A$, we have
\[
\liminf_{n \rightarrow \infty}\frac{A(n)}{n^{\frac{k-2}{k-1}}} > c_{k},
\]
where $c_{k} > 0$ is a constant and $k \ge 3$.
\end{itemize}
\end{conj}

Finally we change the number $2$ in Remark 1 to any integer $k$ and obtain the following result.

\begin{thm}
Let $a_{0} \ge 3$ and $k \ge 4$ be fixed. Let $A = \{a_{0},\ldots\}$ be defined by the following
greedy algorithm: for any integer $n > a_{0}$, $n \in A$ if and only if
 $n \ne kb - a$, where $a$, $b < n$ and $a, b \in A$. Then we have
\[
A = \bigcup_{n=0}^{\infty}[a_{n}, b_{n}],
\]
where $b_{0} = \left \lfloor \frac{ka_{0}}{2} \right \rfloor$, $a_{l} = kb_{l-1} - a_{0} + 1$ and
$b_{l} = \left\lfloor \frac{ka_{l}}{2} \right \rfloor$ for all integers $l\ge 1$.
\end{thm}

\begin{rem}
If one of the conditions $a_0\ge 3$ and $k\ge 4$ does not hold, then some sequences generated by $\{a_{0}\}$ seems to be chaotic,
without nice structure.
\end{rem}

\section{Proof of Theorem 1}

We define the sequence $n_{k}$ recursively.
Let $n_{1} = 1$ and $n_{k+1} =2^{2n_{k}+2}$ for $k=1,2,\ldots$. Define sets
\[
A_{k} = \{n_{k} + 1, n_{k} + 2, \dots{}, n_{k} + 2^{n_{k}+1}\} \cup \{3\cdot 2^{n_{k}}, 4\cdot 2^{n_{k}}, \dots{}, (2^{n_{k}+1}+2)\cdot 2^{n_{k}}\}
\]
and $A = \cup_{k = 1}^{\infty}A_{k}$.

Now we prove that for any integer $n$, there exist $a,b\in A$ with $a<b<n$ such that $a,b,n$ form an arithmetic progression of length three.

Take an integer $k$ such that $n_{k} + 3 \le n < n_{k+1} + 3$.  It is enough to prove that there exist
$a, b\in A_{k}$ with $a < b < n$ such that $a$, $b$, $n$ form an
arithmetic progression of length three.

Case 1. $n_{k} + 3 \le n \le n_{k} + 2^{n_{k}+1}$. In this case, we take $a = n - 2$, $b = n - 1$. Then $a,b\in A$ and $a$, $b$,
$n$ form an arithmetic progression of length three.

Case 2. $n_{k} + 2^{n_{k}+1}+ 1  \le n <n_{k+1}+3$. It follows that $n \le  2^{2n_{k}+2}+2$. Let
\[
c = 2^{n_{k}}\cdot \left\lceil \frac{n}{2^{n_{k}+1}} \right\rceil,
\] where $\lceil x \rceil$ denotes the least integer not less than $x$.
Then
\[
\frac{n}{2} =2^{n_k}\cdot \frac{n}{2^{n_k+1}}\le c<2^{n_k}\cdot \left(\frac{n}{2^{n_k+1}}+1\right)=\frac{n}{2} + 2^{n_{k}}.
\]
Let $d = 2c - n$. Then $0\le d< 2^{n_k+1}$.

Subcase 2.1. $d>n_{k}$. It follows that $d\in A_k$. Noting that $2\le \lceil \frac{n}{2^{n_{k}+1}} \rceil \le 2^{n_k+1}+1$, we have
$$c= 2^{n_{k}}\cdot \left\lceil \frac{n}{2^{n_{k}+1}} \right\rceil\ge 2^{n_k+1}>d$$ and $c\in A_k$.
Take $a = d$, $b = c$. Obviously, $a, b\in A_{k}$, $a < b < n$ and $a$, $b$, $n$ form an arithmetic progression of length three.

Subcase 2.2. $d \le n_{k}$. Let $a = d + 2^{n_{k}+1}$, $b = c + 2^{n_{k}}$. Then
$2b=a+n$. By $d=2c-n\le n_k$ and $n_k+1+2^{n_k+1}\le n$, we have $$b=c+2^{n_k}=\frac{d+n}{2}+2^{n_k}\le \frac{n_k+n}{2}+2^{n_k}< \frac{n}{2}+\frac{n}{2}=n,$$
$$a=2b-n<2b-b=b.$$
Noting that $2^{n_k+1}\le a\le 2^{n_k+1}+n_k$, $b$ is a multiple of $2^{n_k}$ and $3\cdot 2^{n_k}\le b\le (2^{n_{k}+1}+2)\cdot 2^{n_{k}}$, we have $a,b\in A_k$.
Hence, there exist $a, b\in A_{k}$ with
$a < b < n$ such that $a$, $b$, $n$ form an arithmetic progression of length three.

Noting that $\min A_{k+1}=n_{k+1}+1>2^{2n_k+2}$, we have $$A(2^{2n_{k}+2}) \le n_{k} + |A_k|=n_k+2\cdot 2^{n_{k}+1}.$$ Thus
\[
\liminf_{n \rightarrow \infty}\frac{A(n)}{\sqrt{n}} \le \liminf_{k \rightarrow \infty}\frac{A(2^{2n_{k}+2})}{2^{n_{k}+1}} \le \liminf_{k\rightarrow \infty}\frac{n_{k} + 2\cdot 2^{n_{k}+1}}{2^{n_{k}+1}} = 2.
\]
This completes the proof of Theorem 1.

\section{Proof of Theorem 2}

Let
\[
B_{k} = \left\{\sum_{i=0}^{k-1}\varepsilon_{i}4^{i}: \varepsilon_{i}\in \{1, 2\}\right\}.
\]
We will prove that the set
\[
A = \bigcup_{k=1}^{\infty}\left(\bigcup_{i=0}^{8}(i\cdot 4^{k-1} + B_{k})\right)
\]
satisfies the conditions of Theorem 2. Clearly $B_k\subseteq A$ for any integer $k$.

We first prove that for any positive integers $k$ and $n$ with $3\cdot 4^{k-1}\le n<4^k$, there exist integers $a,b\in B_k$ such that $a<b<n$ and $a,b,n$ form an arithmetic progression of length three. Write
\[
n = \sum_{i=0}^{k-1}\mu_{i}4^{i},
\]
where $\mu_{i}\in \{0,1,2,3\}$. Take
\[
a = \sum_{i=0}^{k-1}\varepsilon_{i}^{(1)}4^{i},\quad
b = \sum_{i=0}^{k-1}\varepsilon_{i}^{(2)}4^{i},
\]
where $\varepsilon_{i}^{(1)} = 2$, $\varepsilon_{i}^{(2)} = 1$ if $\mu_{i} = 0$;
$\varepsilon_{i}^{(1)} = 1$, $\varepsilon_{i}^{(2)} = 1$ if $\mu_{i} = 1$;
$\varepsilon_{i}^{(1)} = 2$, $\varepsilon_{i}^{(2)} = 2$ if $\mu_{i} = 2$;
$\varepsilon_{i}^{(1)} = 1$, $\varepsilon_{i}^{(2)} = 2$ if $\mu_{i} = 3$. Since $3\cdot 4^{k-1} \le n < 4^{k}$, it follows that
$\mu_{k-1} = 3$, and so $\varepsilon_{k-1}^{(1)} = 1$,
$\varepsilon_{k-1}^{(2)} = 2$. Hence $a < b < 3\cdot 4^{k-1}\le n$ and $a,b\in B_k$. It is easy to see that $2b=a+n$, and so $a$, $b$, $n$ form an arithmetic progression of length three.

Next we will prove that for any integer $n$, there exist $a,b\in A$ such that $a<b<n$ and $a,b,n$ form an arithmetic progression.

Write
\[
A_{k} = \bigcup_{i=0}^{8}\left(i\cdot 4^{k-1} + B_{k}\right),\quad k=1,2,\ldots.
\]
For any integer $n$, let $3\cdot 4^{t-1}\le n<3\cdot 4^t$ and let
\[
n^{'} =n- \left(\left \lfloor \frac{n}{4^{t-1}}\right \rfloor-3 \right) \cdot 4^{t-1}.
\]
Then we obtain $3\cdot 4^{t-1}\le n' <4^t$. By arguments above, it follows that
there exist integers $a^{'}$, $b^{'} \in B_{t}$ such that $a^{'} <b^{'}<n^{'}$ form an arithmetic progression of length three.
Now let
\[
a = a^{'} + \left(\left \lfloor \frac{n}{4^{t-1}} \right \rfloor - 3 \right) \cdot 4^{t-1},
\]
\[
b = b^{'} + \left(\left \lfloor \frac{n}{4^{t-1}} \right \rfloor - 3 \right) \cdot 4^{t-1}.
\]
Noting that $3\cdot 4^{t-1}\le n<3\cdot 4^t$, we have $0\le \left \lfloor \frac{n}{4^{t-1}} \right \rfloor - 3\le 8$. Hence $a,b\in A_t$ and
$$b= b^{'} + \left(\left \lfloor \frac{n}{4^{t-1}} \right \rfloor - 3 \right) \cdot 4^{t-1}<3\cdot 4^{t-1}+\left(\frac{n}{4^{t-1}}-3\right)\cdot 4^{t-1}=n.$$
By $2b=a+n$, we have $a<b<n$. Therefore, $a,b\in A$, $a<b<n$ and $a$, $b$, $n$ form an arithmetic progression of length three.

In the next step we give an upper estimation of $A(n)$. It is clear that
$$A=\left(\bigcup_{k=1}^{\infty}B_{k}\right)\bigcup \left(\bigcup_{k=1}^{\infty}\left(\bigcup_{i=1}^{8}(i\cdot 4^{k-1} + B_{k})\right)\right).$$
Write
\[
B_1 = \bigcup_{k=1}^{\infty}B_{k},\quad B_2=\bigcup_{k=1}^{\infty}\left(\bigcup_{i=1}^{8}(i\cdot 4^{k-1} + B_{k})\right).
\]
Then $A=B_1\cup B_2$. If

\[
\left(\bigcup_{i=1}^{8}\left(i\cdot 4^{k-1} + B_{k}\right)\right) \cap [1, n] \ne \emptyset,
\]
then we have $4^{k-1} \le n$, and so $k \le \log_{4}n + 1$. It follows that
$$B_2(n)\le \left|\bigcup_{k\le \log_{4}n + 1}\left(\bigcup_{i=1}^{8}(i\cdot 4^{k-1} + B_{k})\right)\right|\le \sum_{k \le \log_{4}n + 1}8|B_{k}| = 8\sum_{k \le \log_{4}n + 1}2^{k}
< 32\sqrt{n}.$$

Let $4^{s-1}\le n<4^s$. Then

\[
B_1(n) \le B_1(4^s-1) = 2^s \le 2^{\log_4 n+1}=2\sqrt{n}.
\]
Hence, we obtain
\[
A(n)\le B_1(n)+B_2(n)  < 34\sqrt{n}.
\]
This completes the proof of Theorem 2.

\section{Proof of Theorem 3}
The proof of Theorem 3 is based on the probabilistic method due to Erd\H os and
R\' enyi. There is an excellent summary of the probabilistic method in the books
\cite{als} and \cite{hal}. Let $P(E)$ denote the probability of an event $E$.
Define the random set $A$ by $$P(n \in A)
= \text{min}\left\{1, c\left(\frac{\log n}{n}\right)^{\frac{1}{k-1}}\right\},$$
where $c$ is a positive constant. Let
\[
\frac{n}{2k} \le u \le \frac{n}{2(k-1)}
\]
be fixed. Let
\[
Y_{n,u} = \{n - iu: 1 \le u \le k - 1\}.
\]
We prove that if $u \ne v$, then $Y_{n,u} \cap Y_{n,v} = \emptyset$.
Otherwise, if $n - iu = n - jv$, then $iu = jv$, where $i \ne j$. We can assume
that $i > j$, thus $$\frac{k-1}{k-2} \le \frac{i}{j} = \frac{u}{v}
\le \frac{k}{k-1},$$ which is impossible. Let $X_{n,u}$ denotes the event
$Y_{n,u} \subset A$. For every $n \ge n_{0}$ we have
\begin{eqnarray*}
&&P(\nexists l: n - l, n - 2l, \dots{}, n - (k - 1)l \in A)\\
&\le &P\left(\bigcap_{\frac{n}{2k} \le u \le \frac{n}{2(k-1)}}\bar{X}_{n,u}\right) =
\prod_{\frac{n}{2k} \le u \le \frac{n}{2(k-1)}}\left(1 - \prod_{i=1}^{k-1}
c\cdot\left(\frac{\log(n-iu)}{n-iu}\right)^{1/(k-1)}\right)\\
&\le& \prod_{\frac{n}{2k} \le u \le \frac{n}{2(k-1)}}\left(1 - \prod_{i=1}^{k-1}
c\cdot\left(\frac{\log n}{n}\right)^{1/(k-1)}\right)
= \prod_{\frac{n}{2k} \le u \le \frac{n}{2(k-1)}}\left(1 - c^{k-1}\frac{\log n}{n}\right)\\
&\le & \prod_{\frac{n}{2k} \le u \le \frac{n}{2(k-1)}}\text{exp}\left(- c^{k-1}
\frac{\log n}{n}\right) \le \text{exp}\left(- c^{k-1}
\frac{\log n}{n}\cdot \frac{n}{2k(k-1)}\right)\\
&\le & \text{exp}(-2\log n) = \frac{1}{n^{2}}.
\end{eqnarray*}
if $c$ is large enough. We will apply the following important lemma.
\begin{lem}\cite[Borel-Cantelli, See p.135]{hal}
Let $E_{1}$, $E_{2}$, ... be a sequence of events in a probability space. If
\[
\sum_{j=1}^{+\infty}P(E_{j}) < +\infty,
\]
\noindent then with probability 1, at most a finite number of the events
$E_{j}$ can occur.
\end{lem}

It follows from Lemma 1 that with probability $1$,
there are only finitely many $n$ such that there does not exist $l$ such that
$n - l, n - 2l, \dots{}, n - (k - 1)l \in A$.
It is easy to see from the method of the proofs of Lemmas 10 and 11 in \cite{hal}, pp. 144 - 145 that
with probability $1$, $A(n) \ll_{k} (\log n)^{1/(k-1)}\cdot n^{\frac{k-2}{k-1}}$.
Thus, with probability $1$, there exist $AP_{k}$ - covering sets $A$ with
$A(n) \ll_{k} (\log n)^{1/(k-1)}\cdot n^{\frac{k-2}{k-1}}$.

\section{Proof of Theorem 4}

Let $I_{l} = [a_{l}, b_{l}]$ and $J_{l} = [b_{l} + 1, a_{l+1} - 1]$.
First we prove that for any $n \in I_{l}$, $a$, $b \in A$ and $a$, $b < n$,
we have $n \ne ka - b$.

Suppose that $n \in I_{l}$ and $n = ka - b$, where $a$, $b \in A$ and $a$, $b < n$. Then if $a \in I_{j}$ for some $j \le l - 1$, we have
\[
kb_{l-1} - a_{0} + 1 = a_{l} \le n = ka - b \le kb_{l-1} - a_{0},
\]
a contradiction. If $a \in I_{l}$, and $b < n$ then
\[
\left \lfloor \frac{ka_{l}}{2} \right \rfloor = b_{l} \ge n = ka - b \ge ka_{l} - (n - 1)
\ge ka_{l} - \left \lfloor \frac{ka_{l}}{2} \right \rfloor + 1 \ge \left \lfloor \frac{ka_{l}}{2} \right \rfloor + 1,
\]
which is also a contradiction.

Hence, for any $n \in I_{l}$, $a$, $b \in A$ and $a$, $b < n$,
we have $n \ne ka - b$.

In the next step, we prove that for any integers $l\ge 0$ and $n \in J_{l}$, there exist $a$, $b \in A$ with $a$, $b < n$ such that $n = ka - b$.

Suppose that $n\in J_{l}$.
For $h=0,1,\ldots,b_l-a_l$, we define
$$J_{l}^{(h)} = k(a_{l} + h) - I_{l} =  \left\{k(a_{l} + h) - i: i\in I_{l}\right\}.$$
It is easy to see that the smallest element of $J_{l}^{(h+1)}$ is
$$\text{min}~J_{l}^{(h+1)} = k(a_{l} + h + 1) - \left\lfloor \frac{ka_{l}}{2} \right\rfloor$$ and the largest element of $J_{l}^{(h)}$ is
$$\text{max}~J_{l}^{(h)} = k(a_{l} + h) - a_{l}.$$ Since $k \ge 4$ and $a_{0} \ge 3$, it follows that for any $h$ with
$0 \le h \le b_{l} - a_{l} - 1$, we have
\[
\text{min}~J_{l}^{(h+1)} - \text{max}~J_{l}^{(h)} = k(a_{l} + h + 1) - \left\lfloor \frac{ka_{l}}{2} \right\rfloor - (k(a_{l} + h) - a_{l}) = k + a_{l} - \left\lfloor \frac{ka_{l}}{2} \right\rfloor \le 1.
\]
It follows that
\[
[b_{l} + 1, kb_l-a_l]\subseteq \bigcup_{h=0}^{b_{l}-a_{l}}J_{l}^{(h)}  .
\]
Hence, for any integer $n\in [b_{l} + 1, kb_l-a_l]$, there exist integers $h$ with $0\le h\le b_l-a_l$ and $i\in I_l$ such that
$n=k(a_l+h)-i$. Clearly $i\le b_l<n$ and $a_l+h \le a_l + (b_l - a_l) = b_l  <n$. Thus we have $i \in A$, $a_{l} + h \in A$. 

It remains to show that for any $kb_{l} - a_{l} + 1 \le n \le kb_{l} - a_{0}$ there exist $a$, $b \in A$, $a$, $b < n$ such that $n = ka - b$.
If $l=0$, then $kb_l-a_l=kb_0-a_0=a_1-1=a_{l+1}-1$, and so $$[b_{l} + 1, kb_l-a_l]=[b_l+1,a_{l+1}-1].$$
Now we suppose that $l\ge 1$.

Let $K_l=\{ka-b:~a\in I_l,~b\in I_0\}$. Since $l\ge 1$, it follows that $a>b$. By $k\ge 4$, we have $a<ka-b$ and $b<ka-b$.
By $k\ge 4$ and $a_0\ge 3$, we have
$$|I_0|=b_0-a_0+1=\left\lfloor \frac{ka_0}{2}\right\rfloor-a_0+1\ge k.$$
It follows that $K_l=[ka_l-b_0,kb_l-a_0]$. By $b_l=\lfloor \frac{ka_l}{2}\rfloor $ and $b_{0} \ge k \ge 4$, we have
$$kb_l-a_l = k\left\lfloor \frac{ka_l}{2}\right\rfloor - a_{l} > k\left (\frac{ka_l}{2} - 1\right) - a_{l} = \frac{k^2}{2}a_l-k-a_l>ka_l-b_0.$$
Hence $$[b_{l} + 1, kb_l-a_l]\cup [ka_l-b_0,kb_l-a_0]=[b_l+1,kb_l-a_o]=[b_l+1,a_{l+1}-1].$$
Therefore, if $n \in J_{l}$, then there exist $a$, $b \in A$ and $a$, $b < n$ such that $n = ka - b$.

This completes the proof of Theorem 4.

\section{Acknowledgement} Part of this work was done during the third author visiting to Budapest
University of Technology and Economics. He would like to thank Dr.
S\'andor Kiss and Dr. Csaba S\'{a}ndor for their warm hospitality.

\end{document}